\title{ ~~\\ A combinatorial identity arising from cobordism theory}
\author{Dion Gijswijt and Pieter Moree}
\documentclass[12pt]{article}
\usepackage{amssymb, latexsym, amsfonts}
\def\@ptsize{2}
\newtheorem{Thm}{Theorem}

\newtheorem{Cor}{Corollary}

\newtheorem{Prop}{Proposition}
\newcommand{\qed}{\hfill $\Box$}

\begin{document}
\date{}
\maketitle
\centerline{\it Dedicated to the memory of Alexander Reznikov}
{\def\thefootnote{}
}
\begin{abstract}
\noindent Let ${\underline \alpha}=(\alpha_1,\alpha_2,\cdots,\alpha_m)\in \mathbb R_{>0}^m$.
Let ${\underline \alpha}_{i,j}$ be the vector obtained from ${\underline \alpha}$
on deleting the entries $\alpha_i$ and $\alpha_j$. Besser and Moree \cite{BM} introduced some invariants and
near invariants  related to the solutions
${\underline \epsilon}\in \{\pm 1\}^{m-2}$ of the linear inequality $|\alpha_i-\alpha_j|<\langle {\underline \epsilon},{\underline \alpha}_{i,j} \rangle  < \alpha_i+\alpha_j,$
where $\langle , \rangle$ denotes the usual inner product and
$\underline  \alpha_{i,j}$ the vector obtained from  $\underline \alpha$ on deleting $\alpha_i$ and $\alpha_j$.
The main result of Besser and Moree \cite{BM} is extended here
to a much more general setting, namely that of certain maps from finite sets to $\{-1,1\}$.
\end{abstract}
\section{Introduction}
Let $m\ge 3$. Let
${\underline \alpha}=(\alpha_1,\alpha_2,\ldots,\alpha_m)\in  \mathbb R_{>0}^m$
and
suppose that there is no
$\underline \epsilon \in  \{\pm 1 \}^m$ satisfying
$\left< \underline \epsilon,\underline
 \alpha \right>=0$.
Let $1\le i<j\le m$. Let $\underline  \alpha_{i,j}\in
\mathbb R_{>0}^{m-2}$ be the vector obtained from  $\underline \alpha$ on deleting $\alpha_i$ and $\alpha_j$. Let
$$
S_{i,j}(\underline \alpha) :=\{\underline \epsilon\in
\{\pm 1    \}^{m-2}: |\alpha_i-\alpha_j|< \left<\underline
 \epsilon,\underline \alpha_{i,j}\right>     <\alpha_i+\alpha_j\}.
$$
Define $N_{i,j}(\underline \alpha) =\sum_{\underline
\epsilon\in S_{i,j}({\underline \alpha})} \prod_{k=1}^{m-2}\epsilon_k$. Theorem 2.1 of \cite{BM} 
states that the reduction of $\# S_{i,j}(\underline \alpha)$ mod $2$ only depends on $\underline \alpha$ and that in case $m$ odd,
$N_{i,j}(\underline \alpha)$ only depends on $\underline \alpha$.
In particular it was shown that for $m\ge 3$ and odd
we have
\begin{equation}
\label{besser}
N_{i,j}(\underline \alpha)=-{1\over 4}\sum_{\underline \epsilon \in \{\pm 1\}^m}{\rm sgn}(\left<\underline
 \epsilon,\underline \alpha\right>)    \prod_{k=1}^m \epsilon_k.
\end{equation}
{}From (\ref{besser}) we of course immediately read off that
if $m\ge 3$ is odd, $N_{i,j}(\underline \alpha)$ does not depend on the choice of $i$ and $j$.\\

\noindent {\bf Example 1.1}. We take ${\underline \beta}_m=(\log
2,\dots,\log p_m)$,
where $p_1,\dots,p_m$ denote the consecutive primes and
put $Q=p_1\cdots p_m$. Then it is not difficult
to show that, for $1\le i<j\le m$,
$$N_{i,j}({\underline \beta}_m)=(-1)^m\sum_{\sqrt{Q/p_i}<n<\sqrt{Q}\atop {\rm gcd}(n,p_ip_j)=1,~P(n)\le p_m}\mu(n),$$
where $P(n)$ denotes the largest prime factor of $n$ and $\mu$ the M\"obius function.
 For $m\ge 2$ put
 $$g(m)={(-1)^{m+1}\over 4}\sum_{d|p_1\cdots p_m}{\rm sgn}({d^2\over
p_1\cdots p_m}-1)
\mu(d),$$
where sgn denotes the sign function.
The fundamental theorem of arithmetic ensures there is no
$\underline \epsilon \in  \{\pm 1 \}^m$ satisfying
$\langle \underline \epsilon,
 {\underline \beta}_m \rangle=0$. By (\ref{besser}) we then infer
 that
if $m\ge 3$ is odd, $N_{i,j}({\underline \beta}_m)=g(m)$ and so does not depend on the choice of $i$ and $j$.
By Remark 2.5 of \cite{BM} we have $g(m)=0$ for $m\ge 2$ and even. The first non-trivial values one finds for $g(m)$ 
are given in the table below.

\begin{center}
\begin{tabular}{|c|c|c|c|c|c|c|c|c|c|c|c|c|c|c|c|c|}
\hline
$m$  & $3$ & $5$ & $7$ & $9$ & $11$ & $13$ & $15$ & $17$ & $19$ & $21$ & $23$\\
\hline
$g(m)$  & $1$ & $-1$ & $3$ & $-8$ & $22$ &
$-53$ & $158$ & $-481$ & $1471$ & $-4621$ & $14612$ \\
\hline
\end{tabular}
\end{center}

\noindent (The value given for $m=15$ corrects the value at p. 471 of \cite{BM}. For
a computer program to evaluate these values see \cite{N}.)\\

\noindent {\bf Example 1.2}. Put $Q(n)=\sum_{d|n,~d\le \sqrt{n}}\mu(d)$. The sequence
$\{Q(0),Q(1),Q(2),\dots\}$ is sequence A068101 of OEIS \cite{Q}.\\ \indent Let $n>1$ be a squarefree integer 
having $k$ distinct prime divisors $q_1,\ldots,q_k$ with $k\ge 2$.
Note that in the previous example we only used that $p_1,\dots,p_m$ are distinct primes. If
we replace them by $q_1,\dots,q_k$ we infer, proceeding as in the previous example, that
$$g_n(k):={(-1)^{k+1}\over 4}\sum_{d|n}{\rm sgn}({d^2\over n}-1)
\mu(d)$$
is an integer that equals zero if $k$ is even. On using that
$\sum_{d|n}\mu(d)=0$ it is seen
that $g_n(k)={(-1)^k\over 2}Q(n)$, whence the following result is inferred: \begin{Prop}
Let $n>1$ be a squarefree number having $k$ distinct prime divisors. Then $$Q(n)=\cases{1 & if $n$ is a prime;\cr
0 &if $k$ is even;\cr
{\rm even} & if $k\ge 3$ is odd.}$$
\end{Prop}
\section{General setup}
We consider a more general quantity $N_{\sigma}(a,b)$ similar
to $N_{i,j}(\underline \alpha)$ so that the latter is a special case of the former.\\
\indent Let $X$ be a finite set. Suppose that
we have a map $\sigma:2^{X}\rightarrow \{-1,1\}$
such that $\sigma(X\backslash A)=\sigma(A)$ for all $A\subseteq X$. We will call such a map $\sigma$ {\it even}.
Let $u,v\in X$ with $u\ne v$.
Define
\begin{equation}
\label{definition}
N_{\sigma}(u,v):=\sum_{A\subseteq X,~u\in A,~v\not\in A\atop
\sigma(A)=\sigma(A+v)}\sigma(A),
\end{equation}
where the summation is over all subsets $A$ of $X$ such
that $u\in A$, $v\not\in A$ and $\sigma(A)=\sigma(A+v)$.
\begin{Thm}
\label{main}
Let $\sigma$ be an even map from $X\rightarrow \{-1,1\}$. Then
$$N_{\sigma}(u,v)={1\over 4}\sum_{A\subseteq X}\sigma(A)$$
and thus in particular
$N_{\sigma}(u,v)$ does not depend on the
 choice of $u$ and $v$.
 \end{Thm}
{\it Proof}. We have
 \begin{eqnarray}
2N_{\sigma}(u,v)
   &=&\sum_{A\subseteq X,~u\in A,~v\not\in A\atop
\sigma(A)=\sigma(A+v)} (\sigma(A)+\sigma(A+v))
   =\sum_{A\subseteq X \atop u\in A, v\not\in A} (\sigma(A)+\sigma(A+v)) \cr &=&\sum_{A\subseteq X\atop u\in A} \sigma(A)
   ={1\over 2} \sum_{A\subseteq X\atop u\in
A}(\sigma(A)+\sigma(X\backslash A)),\cr
   &=&{1\over 2}(\sum_{A\subset X\atop u\in A}\sigma(A)+\sum_{A\subseteq
X\atop
   u\not\in A}\sigma(A))={1\over 2}\sum_{A\subseteq X}\sigma(A),\nonumber
\end{eqnarray}
where we used that there is a bijection between the sets containing $u$ and those not containing $u$, the bijection being taking complementary sets. \qed\\

\noindent {\tt Remark}. In case the cardinality
of $X$ is odd, we can alternatively consider a
map $\tau:2^{X}\rightarrow \{-1,1\}$ such that
$\tau(X\backslash A)=-\tau(A)$ for all $A\subseteq X$.
Then the map $\sigma$ defined by $\sigma(A)=(-1)^{\# A}\tau(A)$ is even and the conditions of Proposition \ref{main} are satisfied.

\section{Examples}
We present three applications of Theorem 1.\\

\noindent {\bf Example 3.1}. Suppose $X=\{x_1,\dots,x_m\}$ and $m\ge 3$. Let $f$ be a map such that
$f(x_j)=\pm 1$ for $1\le j\le m$. Consider the map
$\sigma:2^{X}\rightarrow \{-1,1\}$
defined by
$\sigma(A)=\prod_{a\in A}f(a)$ for $A\subseteq X$. Let us assume that $\prod_{x\in X}f(x)=1$ 
(so that $\sigma$ is an even map). 
Theorem 1 then gives that
$$N_{\sigma}(u,v)=\cases{2^{\# X -2} &if $f(x_j)=1$ for $1\le j\le m$;\cr 0 &otherwise.}$$

\noindent {\bf Example 3.2}. We reprove the main result from \cite{BM} which is reproduced in the present note as (\ref{besser}), where we
now drop the requirement that $\alpha_j>0$ for $1\le j\le m$. Let
$X=\{\alpha_1,\dots,\alpha_m\}$ be a set of cardinality $m$ consisting of real numbers such that there is no
$\underline \epsilon \in  \{\pm 1 \}^m$ satisfying
$\langle \underline \epsilon,\underline
 \alpha \rangle=0$. Let $A$ be any subset of $X$. To $A$ we associate
${\underline \epsilon}=(\epsilon_1,\dots,\epsilon_m)$, where
$\epsilon_j=-1$ if $\alpha_j\in A$ and $\epsilon_j=1$ otherwise.
Let $\sigma(A)={\rm sgn}(\langle \underline \epsilon,\underline
 \alpha \rangle)\epsilon_1\cdots \epsilon_m$. By assumption
 $\langle \underline \epsilon,\underline
 \alpha \rangle\ne 0$ and hence $\sigma(A)\in \{-1,1\}$. Let $i\ne j$. We evaluate $N_{\sigma}(\alpha_i,\alpha_j)$ according to the
 definition (\ref{definition}). We obtain
 that $N_{\sigma}(\alpha_i,\alpha_j)=\sum' {\rm sgn}(
 \langle \underline \epsilon,\underline
 \alpha \rangle)\prod_{k=1}^m \epsilon_k$, where the dash indicates that we sum over those ${\underline \epsilon}\in \{\pm 1\}^m$, where
 $\epsilon_i=-1$, $\epsilon_j=1$ and
 $$-{\rm sgn}(\langle {\underline \epsilon}_{i,j},{\underline
 \alpha}_{i,j} \rangle - \alpha_i+\alpha_j)={\rm sgn}
 (\langle {\underline \epsilon}_{i,j},{\underline
 \alpha}_{i,j} \rangle - \alpha_i-\alpha_j).$$ Note that the latter condition is satisfied iff $\alpha_i-|\alpha_j| <
 \langle {\underline \epsilon}_{i,j},{\underline
 \alpha}_{i,j} \rangle < \alpha_i+|\alpha_j|$. If
  ${\underline \epsilon}\in \{\pm 1\}^m$ satisfies the latter inequality, $\epsilon_i=-1$
  and $\epsilon_j=1$, then $${\rm sgn}(
 \langle \underline \epsilon,\underline
 \alpha \rangle)\prod_{k=1}^m \epsilon_k=-{\rm sgn}(\alpha_j)\prod_{k=1\atop k\ne i,j}^m \epsilon_k.$$
We infer that $$N_{\sigma}(\alpha_i,\alpha_j)=-{\rm sgn}(\alpha_j) \sum_{{\underline \epsilon}\in \{\pm 1\}^{m-2}\atop
\alpha_i-|\alpha_j| < \langle {\underline \epsilon},{\underline \alpha}_{i,j} \rangle < \alpha_i+|\alpha_j|}\prod_{k=1}^{m-2}\epsilon_k.$$
In case $m$ is odd, $\sigma$ is even and Theorem 1 can be applied
(note that $N_{\sigma}(\alpha_i,\alpha_j)=-{\cal N}_{i,j}(\underline \alpha)$) to give the following
corollary.
\begin{Cor}
Let ${\underline \alpha}=(\alpha_1,\alpha_2,\cdots,\alpha_m)\in \mathbb R^m$ and
suppose that there is no
$\underline \epsilon \in  \{\pm 1 \}^m$ satisfying
$\left< \underline \epsilon,\underline
 \alpha \right>=0$. Let $1\le i<j\le m$. Put
$$
{\cal S}_{i,j}(\underline \alpha) :=\{\underline \epsilon\in
\{\pm 1    \}^{m-2}: \alpha_i-|\alpha_j|< \left<\underline
 \epsilon,\underline \alpha_{i,j}\right>     <\alpha_i+|\alpha_j|\}.
$$
Define ${\cal N}_{i,j}(\underline \alpha) ={\rm
sgn}(\alpha_j)\sum_{\underline
\epsilon\in {\cal S}_{i,j}({\underline \alpha})}
\prod_{k=1}^{m-2}\epsilon_k$.
If $m\ge 3$ and $m$ is odd, then
$${\cal N}_{i,j}(\underline \alpha)=-{1\over 4}\sum_{\underline \epsilon \in \{\pm 1\}^m}{\rm sgn}(\left<\underline
 \epsilon,\underline \alpha\right>) \prod_{k=1}^m \epsilon_k=
 h({\underline \alpha}),$$ does not depend on $i$ and $j$. If one of the entries of $\underline \alpha$ is zero, 
 then $h({\underline \alpha})=0$.
 \end{Cor}
In case  ${\underline \alpha}\in \mathbb R_{>0}^m$ it is not immediately clear that this result implies (\ref{besser}). To see that this
is nevertheless true it suffices
to show that under the conditions of Corollary 1 we have
${\cal N}_{i,j}(\underline \alpha)=N_{i,j}(\underline \alpha)$.
If $\alpha_j\le \alpha_i$ this is obvious, so assume that
$\alpha_j> \alpha_i$.
Notice that $\underline \epsilon\in
\{\pm 1    \}^{m-2}$ is in ${\cal S}_{i,j}(\underline \alpha)\backslash {S}_{i,j}(\underline \alpha)$ iff
$\alpha_i-\alpha_j< \left<\underline
 \epsilon,\underline \alpha_{i,j}\right>     <\alpha_j-\alpha_i$. But if $\underline \epsilon$ satisfies 
 the latter inequality, so does $-{\underline \epsilon}$ and both are counted with opposite sign
 in ${\cal N}_{i,j}(\underline \alpha)-{N}_{i,j}(\underline \alpha)$ and consequently
 ${\cal N}_{i,j}(\underline \alpha)={N}_{i,j}(\underline \alpha)$.\\

\noindent {\bf Example 3.3}. Corollary 1 can be
generalised to a higher dimensional
setting.
Instead of numbers $\alpha_1,\dots,\alpha_m$ we can consider
points ${\underline \alpha_1},\dots,{\underline \alpha_m}$
with $\underline \alpha_i\in \mathbb R^n$ and $n\ge 2$.
We assume that $\pm {\underline \alpha_1}\pm \cdots \pm {\underline \alpha_m}\ne {\underline 0}$.
Let us define $B$ to be the matrix
with ${\underline \alpha_j}$ as $j$th row for $1\le j\le m$.
Choose a hyperplane $H$ through the origin not containing any of
the points $\pm {\underline \alpha_1}\pm \cdots \pm {\underline \alpha_m}$ (the assumption that $\pm {\underline \alpha_1}\pm \cdots \pm {\underline \alpha_m}\ne {\underline 0}$ ensures
that this is possible).
Let ${\underline n}\not\in H$ be on the normal of this hyperplane. Let $A$ be any subset of $X$. To $A$ we associate
${\underline \epsilon}=(\epsilon_1,\dots,\epsilon_m)$, where
$\epsilon_j=-1$ if ${\underline \alpha}_j\in A$ and $\epsilon_j=1$ otherwise.
Let $\sigma(A)={\rm sgn}(\langle \underline n, \underline \epsilon B \rangle)\epsilon_1\cdots \epsilon_m$. The assumption on $H$ implies that
 $\langle \underline n,\underline \epsilon B \rangle\ne 0$ and hence
$\sigma(A)\in \{-1,1\}$.
Choose two points ${\underline \alpha_i}$ and ${\underline \alpha_j}$, $i\ne j$.
Let $V$ be the hyperplane with normal ${\underline n}$ containing
${\underline \alpha_i}-{\underline \alpha_j}$ and $W$ be the hyperplane with normal ${\underline n}$ containing ${\underline \alpha_i}+{\underline \alpha_j}$.
We define the weight $w(\underline \alpha)$ of
a point ${\underline \alpha}$ of
the form ${\underline \alpha}=\sum_{1\le k\le m\atop k\ne i,~k\ne
j}\epsilon_k
{\underline \alpha_k}$ with ${\underline \epsilon}_{i,j}\in \{\pm 1\}^{m-2}$ to be $\prod_{1\le k\le m\atop k\ne i,~k\ne j}
\epsilon_k$. Note that our choice of ${\underline n}$ ensures that none of these points is in $V$ or $W$. Then let $M(i,j)$ be the sum of the weights of all points $\sum_{1\le k\le m\atop k\ne i,~k\ne j}\epsilon_k
{\underline \alpha_k}$ that
are in between $V$ and $W$ and
for which ${\underline \epsilon}_{i,j}\in \{\pm 1\}^{m-2}$. If $m\ge 3$ is odd, then $\sigma$ is an even map.
It is not difficult to show that $N_\sigma({\underline \alpha}_i,{\underline \alpha}_j)=\pm M(i,j)$, where
the sign is independent of $i$ and $j$.
Theorem 1 applies and we infer that $M(i,j)$ is independent of the choice of $i$ and $j$.\\

\noindent {\bf Acknowledgement}. We thank Tony Noe for pointing out a typo in \cite{BM} and for 
providing us with the table given in
this note.\\ 
\indent This note has its source in a question posed
by the late Alexander Reznikov to Amnon Besser and the 
second author in the summer of 1997, whilst all three of them were enjoying
the hospitality of the MPI in Bonn. Reznikov came to this question on the basis of computations (together with Luca Migliorini) in the cobordism theory of the moduli space of
polygons. The second author remembers Alexander Reznikov as a
very original and creative mathematician and an intriguing and interesting personality.\\
\indent The research of the second author was made possible thanks to Prof. E. Opdam's PIONIER-grant from the 
Netherlands Organization for Scientific Research (NWO).\\

\medskip\noindent {\footnotesize Dion Gijswijt\\ Korteweg-de Vries Institute\\
Plantage Muidergracht 24\\ 1018 TV Amsterdam\\ 
The Netherlands\\ e-mail: {\tt gijswijt@science.uva.nl}\\

\noindent Pieter Moree\\ Max-Planck-Institut f\"ur Mathematik\\ Vivatsgasse 7\\ D-53111 Bonn\\ Germany\\
e-mail: {\tt moree@mpim-bonn.mpg.de}}

\end{document}